\documentclass[12pt,leqno]{amsart}
\usepackage{amscd,times,fullpage}
\usepackage[utf8]{inputenc}
\usepackage{amsmath}
\usepackage{amsthm}
\usepackage{mathtools}
\usepackage{amsfonts}
\usepackage{comment}
\usepackage{enumerate}
\usepackage{url}
\usepackage{tikz}
\usetikzlibrary{cd}
\usepackage{hyperref}

\hypersetup {
	pdfpagemode = {UseNone},
	pdftitle = {Maximally non-integrable almost complex structures: an h-principle and cohomological properties},
	pdfauthor = {Rui Coelho, Giovanni Placini, Jonas Stelzig},
	pdflang = {en-UK}
} 

\newtheorem{thmintro}{Theorem}
\newtheorem{corintro}{Corollary}[thmintro]

\newtheorem{prop}{Proposition}
\newtheorem{lem}[prop]{Lemma}
\newtheorem{defin}[prop]{Definition}
\newtheorem{definition}[prop]{Definition}
\newtheorem{thm}[prop]{Theorem}
\newtheorem{cor}[prop]{Corollary}
\theoremstyle{remark}
\newtheorem{remark}[prop]{Remark}
\newtheorem{rem}[prop]{Remark}
\newtheorem{ex}[prop]{Example}
\newtheorem{obs}[prop]{Observation}

\newtheoremstyle{colon}%
{}
{}
{\itshape}
{}
{\bfseries}
{:}
{ }
{}

\theoremstyle{colon}

\title[Maximally non-integrable almost complex structures]{Maximally non-integrable almost complex structures: \\ an
	$h$-principle and cohomological properties}
\author{R.~Coelho}
\address{}
\email{}

\author{G.~Placini}
\address{Dipartimento di Matematica e Informatica, Universit\`a di Cagliari, Via~Ospedale~72, 09124 Cagliari, Italy}
\email{giovanni.placini@unica.it}

\author{J.~Stelzig}
\address{Mathematisches Institut, Ludwig-Maximilians-Universit\"at
	M\"unchen, Theresienstr.~39, 80333 M\"unchen, Germany}
\email{jonas.stelzig@math.lmu.de}


\newcommand{\C}{\mathbb{C}}
\newcommand{\HH}{\mathbb{H}}

\newcommand{\R}{\mathbb{R}}
\newcommand{\Z}{\mathbb{Z}}
\newcommand{\mubar}{{\bar\mu}}
\newcommand{\lra}{\longrightarrow}
\newcommand{\fg}{\mathfrak{g}}

\newcommand{\cA}{\mathcal{A}}

\newcommand{\rk}{\operatorname{rk}}
\newcommand{\Hom}{\operatorname{Hom}}

\newcommand{\pr}{\operatorname{pr}}
\newcommand{\im}{\operatorname{im}}
\newcommand{\rank}{\operatorname{rank}}
\newcommand{\coker}{\operatorname{coker}}
\newcommand{\del}{\partial}
\newcommand{\de}{\mbox{d}}
\newcommand{\delbar}{{\bar\partial}}
\newcommand{\lie}{\mathcal{L}}
\newcommand{\cH}{\mathcal{H}}

\newcommand{\J}{\mathcal{J}}
\newcommand{\ZJ}{Z^{(1)}}
\newcommand{\NN}{\mathcal{N}}

\begin{document}

	\begin{abstract}
		We study almost complex structures with lower bounds on the rank of the Nijenhuis tensor. Namely, we show that they satisfy an $h$-principle. As a consequence, all parallelizable manifolds and all manifolds of dimension $2n\geq 10$ (respectively $\geq 6$) admit a almost complex structure whose Nijenhuis tensor has maximal rank everywhere (resp. is nowhere trivial). For closed $4$-manifolds, the existence of such structures is characterized in terms of topological invariants. Moreover, we show that the Dolbeault cohomology of non-integrable almost complex structures is often infinite dimensional (even on compact manifolds).
	\end{abstract}
	\maketitle
	\section{Introduction}
	A basic local invariant of an almost complex structure $J$ on a manifold $X$ is its Nijenhuis tensor or, equivalently, the $(-1,2)$-component of the exterior differential:
	\[
	N_J:=\pr^{0,2}\circ~ d\colon A^{1,0}_X\lra A^{0,2}_X.
	\]
	Therefore almost complex manifolds can be stratified by the pointwise rank of $N_J$. By the celebrated theorem of Newlander-Nirenberg, $J$ is integrable if and only if $N_J\equiv 0$. Manifolds equipped with integrable almost complex structures constitute an extremely interesting class and their properties are largely studied. It seems then natural to ask in which terms one can describe the other extremes or whether they satisfy similar properties. 
	
	Motivated by these questions in this paper we begin to investigate the two arguably most natural cases. Namely, we focus on the spaces of almost complex structures $J$ for which $N_J$ has maximal rank at every point and for which $N_J$ is non-trivial at every point.  Almost complex structures of the former type are said to be maximally non-integrable while those of the latter type are called everywhere non-integrable. Prominent examples of maximally non-integrable almost complex structures arise from (strictly) nearly Kähler six manifolds (like $S^6$) \cite{bryant06,verbitsky08, verbitsky11} and Twistor spaces \cite{cahengerardgutthayyani20,cahenguttrawnsley21}. In \cite{muskarov86}, it was proven that in dimension $\geq 6$ maximally non-integrable almost complex manifolds do not admit any local holomorphic functions. Constant (and in particular maximal) rank conditions on $N_J$ also appear naturally in the study of the recently introduced Dolbeault cohomology of almost complex manifolds \cite{ciriciwilson18}. For example, the Frölicher spectral sequence degenerates at $E_2$ for maximally non-integrable almost complex $4$- and $6$-manifolds. In the present article, we will focus on two aspects: existence and cohomological properties.
	
	\subsection{Existence} 
	Since a generic matrix has maximal rank, one may be tempted to think that any almost complex structure can be deformed into a maximally non-integrable one. However, the existence of a full-rank bundle map from $A^{1,0}_X$ to $A_X^{0,2}$ gives a necessary condition for a homotopy class of almost complex structures to contain maximally non-integrable representatives. Using Gromov's $h$-principle for ample differential relations, we show that this is indeed the only obstruction. In fact, we prove something more. Let us write $\cH$ for the space of pairs $(J,N)$ where $J$ is an almost complex structure on $X$ and $N$ is a global section of the bundle $\Hom_\C(A^{1,0}_X,A^{0,2}_X)$ (where the bigrading is taken according to $J$). Then the first result of this paper is the following:
	\begin{thmintro}\label{thm: h-pple}
		For any manifold $X$ and any integer $r$, sending an almost complex structure $J$ to its Nijenhuis tensor $N_J$ induces a homotopy equivalence
		\[
		\{J \mid \rank_\C (N_J)_x\geq r~\forall x\in X\}\lra \{(J,N)\in \cH\mid \rank_\C (N_x) \geq r~\forall x\in X\}.
		\]
	\end{thmintro}
	Therefore, if the dimension of the manifold $X$ is sufficiently large, the existence of maximally non-integrable or everywhere non-integrable structures in a given homotopy class is not obstructed. In fact a transversality argument yields:
	\begin{corintro}\label{corintro: novacs everywhere}
		In dimension $2n\geq 10$ (resp. $\geq 6$) any almost complex structure on a $2n$-dimensional manifold is homotopic to a maximally non-integrable one (resp. to an everywhere non-integrable one).
	\end{corintro}
	The existence of maximally non-integrable structures is not guaranteed in lower dimensions. In fact, it imposes non-trivial topological conditions on the manifold.
	In particular, in dimension $4$ the Nijenhuis tensor has complex rank $0$ or $1$ at every point, thus everywhere non-integrability and maximal non-integrability are equivalent conditions. In this case it turns out that the existence of such structures is a purely topological question:
	\begin{corintro}\label{corintro: dim 4}
		On a compact oriented $4$-manifold $X$, the following conditions are equivalent:
		\begin{enumerate}[(a)]
			\item $X$ admits a maximally non-integrable almost complex structure inducing its orientation.\label{Eq4DimMaNIACS1}
			\item $X$ admits an almost complex structure inducing its orientation and satisfies $5\chi+6\sigma=0$.\label{Eq4DimMaNIACS2}
			\item $X$ satisfies $5\chi+6\sigma=0$,  $\sigma\equiv0\pmod4$ and if its intersection form is definite, then $b_+=0$.\label{Eq4DimMaNIACS3}
		\end{enumerate}
		Moreover, if the above conditions are satisfied, any almost complex structure inducing the given orientation is homotopic to a maximally non-integrable almost complex structure.
	\end{corintro}
	
	Also in dimension $6$ there is a topological characterization of the existence of maximally non-integrable complex structures:
	\begin{corintro}\label{corintro: dim 6}
		Let $X$ be an oriented 6-manifold.
		\begin{enumerate}
			\item After possibly passing to a finite cover, any almost complex structure on $X$ with $3c_1=0$ (or more generally $c_1$ torsion) is homotopic to a maximally non-integrable one.
			\item $X$ admits a maximally non-integrable almost complex structure if and only if it admits a $\mathrm{spin}$-structure.
		\end{enumerate}
	\end{corintro}
	The second statement should be compared to the result of Wall \cite[Thm 9]{wall66} that $X$ admits an almost complex structure if and only if it admits a $\mathrm{spin}^{c}$-structure.

	Regardless of the dimension, we find that parallelizable manifolds of even dimension always admit maximally non-integrable almost complex structures. In fact:
	\begin{corintro}\label{corintro: parallelizable}
		On a parallelizable manifold of dimension $2n$, any almost complex structure for which the tangent bundle is trivial as a complex vector bundle is homotopic to a maximally non-integrable almost complex structure.
	\end{corintro}
	This corollary applies in particular to left-invariant almost complex structures on homogeneous manifolds. We complement these results with several examples in the main body of the text.

	\subsection{Cohomological properties}
	Cirici and Wilson \cite{ciriciwilson18} proposed a generalization of the Fr\"olicher spectral sequence to the setting of almost complex manifolds. 
	In particular, they defined bigraded vector spaces $H_{Dol}^{\bullet,\bullet}(X)$ generalizing Dolbeault cohomology. 
	This generalization satisfies some desirable properties. 
	Namely, for $X$ compact, the pages of the Fr\"olicher spectral sequence have to become finite-dimensional at some stage because it converges to de Rham cohomology. 
	Moreover, the spaces $H^{p,0}_{Dol}(X)$ and $H^{p,n}_{Dol}(X)$ are known to be finite-dimensional for all $p$, where $X$ is $2n$-dimensional.
	The analogous result is known to hold for all bidegrees in the integrable setting. In view of this, it is a natural question whether finite dimensionality holds for all $H^{p,q}_{Dol}(X)$.  
	We answer this in the negative and show that the Dolbeault cohomology will very often be infinite dimensional:
	
	\begin{thmintro}\label{thm: mani => HDol infty}
		Let $(X,J)$ be an almost complex $4$ or $6$-manifold (not necessarily compact) such that $J$ is maximally non-integrable at a point $x\in X$. Then $H_{Dol}(X)$ is infinite dimensional.
	\end{thmintro}
	In particular, this Theorem applies to $S^6$ with the almost complex structure coming from octonian multiplication and any other nearly Kähler $6$-manifold. By taking products, one obtains examples of almost complex manifolds of arbitrary dimension $2n\geq 4$ with infinite dimensional Dolbeault cohomology.
	Since on an almost complex $4$-manifold $N_J$ has rank $0$ or $1$, the result above specializes to the $4$-dimensional case as:
	\begin{corintro}
		Let $(X,J)$ be an almost complex compact $4$-manifold. Then:
		\[
		\dim H_{Dol}(X,J)<\infty\Longleftrightarrow J\text{ is integrable.}
		\]
	\end{corintro}
	The next two Corollaries were pointed out to us by Joana Cirici.
	As mentioned before, in \cite{ciriciwilson18}, it was shown that all differentials in the Fr\"olicher spectral sequence of maximally non-integrable structures on $4$- and $6$-manifolds vanish from page $2$ onwards. Since de Rham cohomology is finite dimensional for compact manifolds, one may complement this as follows.
	\begin{corintro}
		The Fr\"olicher spectral sequence of a compact $4$ or $6$-dimensional maximally non-integrable almost complex manifold degenerates exactly at page $2$ (and not at page $1$).  
	\end{corintro}
	Note that there are examples of nilmanifolds with maximally non-integrable left-invariant almost complex structures for which the left-invariant Fr\"olicher spectral sequence does degenerate at $E_1$, see \cite{ciriciwilson18}. Further, recall that for nilmanifolds with a left-invariant (integrable) complex structure, it is conjectured that the Dolbeault cohomology can be computed using left-invariant forms. Because left-invariant cohomology is finite-dimensional, the analogous statement in the non-integrable case fails drastically.
	
	\begin{corintro}
		On $4$- or $6$-dimensional nilmanifolds with left-invariant maximally non-integrable almost complex structures, the Dolbeault cohomology can never be computed using left-invariant forms only.
	\end{corintro}
	Finally, recall that on any complex manifold, there exist the Bott-Chern and Aeppli cohomology groups, which fit into the diagram
	\[
	\begin{tikzcd}
		&H_{BC}(X)\ar[ld]\ar[d]\ar[rd]&\\
		H_{\delbar}(X)\ar[rd]\arrow[Rightarrow]{r}&H_{dR}(X)\ar[d]&H_{\del}(X)\ar[ld]\arrow[Rightarrow]{l}\\
		&H_A(X)
	\end{tikzcd}
	\]
	where horizontal arrows represent spectral sequences and all maps are induced by the identity on forms.
	The work \cite{ciriciwilson18} provides a generalization of the horizontal strip of the diagram to the almost complex setting. In Appendix A, we define metric-independent generalizations for $H_{BC}$ and $H_A$ of general almost complex manifolds that maintain many desirable properties from the integrable setting, complementing \cite{ciriciwilson18} and recent harmonic definitions in \cite{tardinitomassini20}. Like Dolbeault cohomology, they can be infinite dimensional on compact manifolds.\\
	
	\subsection{Structure of the paper}
	In Section~\ref{SecExistence} we discuss existence of maximally non-integrable and everywhere non-integrable almost complex structures. In particular, in Section~\ref{SubSecTopology} we introduce the notation and draw some necessary topological conditions for low dimensional manifolds. The proof of Theorem~\ref{thm: h-pple} is carried out in Section~\ref{SecThmA} while Section~\ref{SecApplications} contains the proofs of its corollaries and further consequences. In the last part of Section~\ref{SecExistence} we present some explicit constructions of maximally non-integrable and, more generally, constant rank almost complex structures.
	
	Section~\ref{SecDolbeault} focuses on the cohomological properties of maximally non-integrable almost complex structures. Namely, it contains the proof of Theorem~\ref{thm: mani => HDol infty} as well as explicit examples of this instance. 
	
	In the last part of the paper, Appendix~\ref{AppendixBottChern}, we define metric-independent generalizations for Bott-Chern and Aeppli cohomology. Moreover, we discuss some of their properties and infinite dimensionality.
	
	\bigskip
	
	\noindent\textbf{Acknowledgements:} We would like to thank Michael Albanese, Joana Cirici, Matthias Kemper, Dieter Kotschick, Aleksandar Milivojevic, Jean Ruppenthal, Nicoletta Tardini, Thomas Vogel and Scott Wilson for many discussions and remarks.
	
	\section{Existence}\label{SecExistence}
	\subsection{Maximally non-integrable almost complex manifolds}\label{SubSecTopology}  
	We begin by recalling some known facts on almost complex structures for future reference and to set the notation. Let $X$ be a manifold of dimension $2n$ and suppose $J$ is an almost complex structure on $X$. Denote by $T^{1,0}X$ and $T^{0,1}X$ the $i$ and $-i$ eigenspaces of $J$  in the complexified tangent bundle $TX\otimes \C$. The dual bundles will be denoted by $A_X^{1,0}$, resp. $A_X^{0,1}$ and we write $A^{p,q}_X:=\Lambda^p A^{1,0}_X\otimes \Lambda^qA_X^{0,1}$ for the bundle of $(p,q)$-forms. The space of sections will be denoted by $\cA^{p,q}_X:=\Gamma(X,A_X^{p,q})$. The Nijenhuis tensor of $J$ can be defined both as the tensor $N_J\colon TX\otimes TX\lra TX$ given by 
	\begin{equation}\label{EqNijTensVec}
		N_J(X,Y)=[X,Y]-[JX,JY]+J\left([JX,Y]+[X,JY]\right)
	\end{equation}
	and as the map $\mubar_J\colon A^{1,0}_X\lra A^{0,2}_X$ given by $(-1,2)$-component of the exterior derivative, that is,
	\begin{equation}\label{EqNijTensVForm}
		\mubar_J(\alpha)=\mbox{pr}^{0,2}\circ\,\de(\alpha),
	\end{equation}
	where $\mbox{pr}^{0,2}\colon A^2_X\lra A^{0,2}_X$. If we restrict the argument of (the complexification of) $N_J$ to $(0,1)$-vectors, one may consider $N_J$ as a global section of $A_X^{0,2}\otimes T^{1,0}X=\Hom(A_X^{1,0},A_X^{0,2})$. Under this identification, $N_J= 4\mubar_J$. We will switch freely between both descriptions.
	
	As explained in the introduction we focus on almost complex structures whose Nijenhuis tensor satisfies certain properties.
	We can now define the main object of our interest.
	\begin{defin}\label{DefMaNIACS}
		An  almost complex structure $J$ on a manifold $X$ is called maximally non-integrable if its Nijenhuis tensor $N_J$ or, equivalently, $\mubar_J$ has maximal rank at all $x\in X$.
	\end{defin}
	Moreover, we will be interested in another class of almost complex structures.
	\begin{defin}\label{Defeniacs}
		An almost complex structure $J$ on a manifold $X$ is called everywhere non-integrable if its Nijenhuis tensor $N_J$ or, equivalently, $\mubar_J$ is non-trivial at all $x\in X$.
	\end{defin}
	
	A simple dimension count yields $\rank_\C A^{1,0}_X=n$ and $\rank_\C A^{0,2}_X=\binom{n}{2}$. Therefore we have the following cases:
	\begin{equation}
		J\mbox{ is maximally non-integrable}\Longleftrightarrow
		\begin{cases}
			\mubar_J\mbox{ is surjective} &\mbox{if } n=2\\
			\mubar_J\mbox{ is an isomorphism} &\mbox{if } n=3\\
			\mubar_J\mbox{ is injective} &\mbox{if } n\geq4
		\end{cases}
	\end{equation}
	
	For maximally non-integrable structures, we get an exact sequence of vector bundles
	\[
	0\lra \ker\mubar\lra A^{1,0}_X\lra A^{0,2}_X\lra \coker\mubar\lra 0
	\]
	which yields the following obstructions on the Chern classes of $J$:
	
	\begin{cor}\label{CorObstructions}
		Let $(X,J)$ be a maximally non-integrable almost complex manifold of dimension $2n$. Then the following holds:
		\begin{equation*}
			\begin{cases}
				5 \chi+ 6 \sigma=0&\mbox{if }n=2\\
				3c_1=c_1^2=c_1c_2=0&\mbox{if }n=3\\
				c_3 + 8c_1^3 + c_1c_2=c_4 - c_1c_3 + c_1^2c_2=0&\mbox{if }n=4
			\end{cases}
		\end{equation*}
		
	\end{cor}
	The cases $n=2$ and $n=3$ are due to Armstrong \cite{armstrong97} and Bryant \cite{bryant06} respectively. Notice that the necessary conditions for $8$-manifolds impose topological constraints. In fact, substituting $c_1c_2=-c_3-8c_1^3$ in the last condition we get 
	$$c_4=2c_1c_3+8c_1^4\, .$$
	Therefore the Euler characteristic $\chi$ of a maximally non-integrable almost complex $8$-manifold must be even. Moreover, substituting the same equality in the signature formula
	\[
	\sigma=\frac{1}{45}(3c^2_2-14c_1c_3+14c_4-c_1^4+4c_1^2c_2)
	\]
	yields $3|\chi$. We conclude that the Euler characteristic of a maximally non-integrable almost complex $8$-manifold is divisible by $6$. For $n\geq 5$, the rank of the cokernel of $\mubar$ becomes very large and one does not obtain any constraints on the Chern classes.
	
	\subsection{Proof of Theorem~\ref{thm: h-pple}}\label{SecThmA}
	\subsubsection{The formal Nijenhuis tensor}\label{NijenhuisTensors}
	
	Let $X$ be a $2n$-dimensional smooth manifold.
	We denote by
	\[Z:=\{J\in \mathrm{End}(TX)\,\vert\,J^2=-1\}\overset{\pi}{\lra} X,\]
	the (metric independent) twistor space of $X$. Over $Z$, we have the bundles $V:=\ker(D\pi)$ and $T:=\pi^{\ast}(TX)\xrightarrow{p} Z$. 
	Both come equipped with tautological complex structures $\tilde\J$ and $\J$, respectively, defined as follows.
	Given a point $J_x\in Z$, $x\in X$, the fiber of $V$ over $J_x$ can be described as \[
	V_{J_x}=\{A\in \mathrm{End}(T_xX)\,\vert\,J_xA+AJ_x=0\}
	\]
	and we set
	\[
	\tilde{\J}|_{V_{J_x}}(A):=J_x\circ A,
	\] 
	whereas $T_{J_x}\cong T_xX$ and $\J$ is given by
	\[
	\J|_{T_{J_x}}:= J_x.
	\]
	With the complex structure $\J$ comes a splitting of the complexification $T_\C$ of $T$ into complex bundles given as usual by the $\pm i$ eigenvalues of $\J$:
	\[
	T_\C=T^{1,0}\oplus T^{0,1}.
	\]
	Writing $A^{p,q}$ for the duals of the bundles $\Lambda^{p}T^{1,0}\otimes \Lambda^{q}T^{0,1}$, we now define
	\begin{equation}\label{EqNijenhuisBundle}
		H:=A^{0,2}\otimes_\C T^{1,0},  
	\end{equation}
	the bundle of \emph{linear Nijenhuis tensors}.
	Note that for a given local section $J:U\subset X\rightarrow Z$ of $Z$, the splitting of $T_\C U$ induced by $J$ is compatible with that of $T$ in the sense that there is a natural bundle embedding of $T_\C U$ into the restriction of $T_\C$ to the image of $J$ in $Z$.
	Moreover, the Nijenhuis tensor $N_J$ of $J$ is naturally a local section of the pull-back bundle $J^{^\ast}(H)$.
	
	Denote by $\ZJ$ the space of $1$-jets of the fibration $\pi:Z\rightarrow X$. 
	Recall that $\pi^{(1)}:\ZJ\to Z$ has the structure of an affine bundle.
	We will now define a bundle map $\NN:\ZJ\rightarrow H$.
	Fix $x\in X$ and $J_x\in Z_x$.
	A point in $\ZJ$ over $J_x$ consists of the 1-jet $j$ of an almost complex structure whose value at $x$ equals $J_x$.
	Picking an extension $J$ for $J_x$ with $1$-jet $j$, put $\NN(j) = N_J$, the Nijenhuis tensor of $J$.
	Since the Nijenhuis tensor of $J$ depends only on $J$ up to first order, $\NN$ is a well-defined bundle map. This setup is summarized in the following diagram:
	\[
	\begin{tikzcd}[column sep=0em]\label{twistor jets diagram}
		& \ZJ \ar[rrrr,"\NN"] \ar[drr,"\pi^{(1)}",swap] &&&& H \ar[lld,"p_H"]\\
		& &&Z\ar[d,"\pi"]&&\\
		& &&X&& 
	\end{tikzcd}
	\]
	Let us describe $\NN$ restricted to a fixed fiber over $J_x\in Z$ more explicitly. We note that the Lie-bracket induces a map
	\begin{equation}\label{eqn: bracket}
		TX^{(1)}\times TX^{(1)}\lra TX
	\end{equation}
	and the action of endomorphisms on the tangent bundle induces a map
	\begin{equation}\label{eqn: endomorphisms}
		\mathrm{End}_\R(TX)^{(1)}\times TX^{(1)}\lra TX^{(1)}
	\end{equation}
	If we only care about $\NN$ at a fiber, we may assume without loss of generality that $X$ is a small ball, s.t. $TX$ is trivial. Recall that as soon as one has a trivial fibration $\pr:V:=X\times F\to X$, the fiber of the $1$-jet space $V^{(1)}\to X$ over $x\in X$ may be canonically identified with the vector bundle $\Hom_\R(\pr^{-1}T_xX, TF)$. Using this identification, we write elements $v\in TX^{(1)}_x$ as $(v_0,v_1)$ with $v_0\in T_xX$ and $v_1=v_1(\_)\in \Hom_\R(T_xX, T_xX)$. Similarly, we write elements in $A\in \mathrm{End}_\R(TX)^{(1)}_x$ as $A=(A_0,A_1)$ with $A_0\in \mathrm{End}_\R(T_xX)$ for some $x\in X$ and $A_1:=A_1(\_)\in \Hom_\R(T_xX,\mathrm{End}(T_xX))$ (in particular, this applies to $A\in Z^{(1)}$). With these identifications, we may write \eqref{eqn: bracket} and \eqref{eqn: endomorphisms} as
	\begin{equation}
		[v,w]=w_1(v_0)-v_0(w_1)
	\end{equation}
	and
	\begin{equation}
		A(v)= (A_0(v_0), A_1(\_)(v_0) + A_0\circ v_1(\_))
	\end{equation}
	Using these formulae, one verifies that for any $(J_x,\lambda)\in Z^{(1)}_{J_x}$ where $\lambda\in \Hom_\R(T_xX,V_{J_x})$, and $v,w\in T^{0,1}_{J_x}$ one has
	\[
	\NN(J_x,\lambda)(v,w)=i(\lambda(v)(w)-\lambda(w)(v)) + J_x(\lambda(v)(w)-\lambda(w)(v)),
	\]
	where $\lambda$ is extended $\C$-linearly. The space $V_{J_x}$ is equipped with the complex structure $\tilde{\J}$ introduced above and therefore so is $Z^{(1)}_{J_x}=\{J_x\}\times\Hom_\R(T_xX,V_{J_x})$, by post-composition. With respect to this complex structure, $\NN$ is complex linear. In fact:
	\begin{align*}
		i \NN(J_x,\lambda)(v,w)&=(-1)(\lambda(v)(w)+i\lambda(w)(v))+i(J_x\lambda(v)(w)-J_x\lambda(w)(v))\\
		&=\NN(j,\tilde{\J} \lambda)(v,w)
	\end{align*}
	we conclude:
	\begin{lem}
		The restriction of $\NN$ to the fibers of $\pi^{(1)}$ is a holomorphic map.
	\end{lem}
	\begin{rem}
		The vertical tangent bundle $\ker D\pi^{(1)}$ can be identified with $(\pi^{(1)})^{*}\Hom(\pi^*TX,V)$. The complex structure it inherits from $\tilde{J}$ by post-composition coincides with the one described above on the fibers of the bundle $\pi^{(1)}$ (after identifiying the fibers, which are affine spaces, with their tangent space at any point). I.e., the fiberwise complex structure is independent of the choice of a local trivialization.
	\end{rem}
	\subsubsection{The h-principle}
	On $H$ we define the sub-bundles
	\[
	H^{\geq i}=\{\lambda\in H\,|\, \mathrm{rank}(\lambda)\geq i\}\subseteq H.
	\]
	We are interested in the (open) partial differential relations  
	\[
	\mathcal{R}^{\geq i} := \NN^{-1}(H^{\geq i})\subseteq Z^{(1)}.
	\]
	We will show that these partial differential relations are \emph{ample}, which will imply the result according to \cite{gromov86}, see also \cite{eliashbergmishachev02}. Let $F:=Z^{(1)}_{J_x}$ be any fiber of $\pi^{(1)}$. Recall that to show ampleness of a partial differential relation which sits in a space of $1$-jets, it is enough to show that the complement of its intersection with any principal subspace is a stratified set of even (real) codimension. We refer to \cite{eliashbergmishachev02} for details and definitions and only spell out what these notions mean in our context:
	
	One may identify $F$ with the vector space $\Hom_\R(T_xX,V_{J_x})$. Given any hyperplane $\tau\in T_xX$ and a linear map $l:\tau\to V_{J_x}$ the affine subspace 
	\[
	F_\tau^{l}:=\{L\in F\mid L|_{\tau}=l\}\subseteq F
	\]
	is called a principal subspace of $F$. The tangent space to $F$ at any point is again $\Hom_\R(T_xX, V_{J_x})$ and the tangent space to $F_\tau^{l}$ is $F_\tau^0$, which is $\tilde\J$-stable. In particular, $F_\tau^l$ is a complex submanifold of $F$. Therefore, $\NN$ restricts to a holomorphic map on $F_\tau^l$ and $\NN^{-1}(H\setminus \tilde{H})\cap F_\tau^l$ is a complex analytic subset of $F$, so it has even (real) codimension. Denoting by $\Gamma(X,Z)^{\geq i}$ the set of sections $J:X\to Z$ with $\rk N_J\geq i$ and by $\Gamma(X,Z^{(1)})^{\geq i}$ the set of sections $J:X\to Z^{(1)}$ with  $\rk N_J\geq i$, Gromov's h-principle for ample differential relations yields a homotopy equivalence
	\begin{align*}
		\Gamma(X,Z)^{\geq i}\lra \Gamma(X,Z^{(1)})^{\geq i}.
	\end{align*}
	To complete the proof of Theorem~\ref{thm: h-pple}, we note that $\NN$ is surjective by the following Lemma:
	\begin{lem}[{\cite[Theorem~5]{kruglikov98}}]\label{lem: Nijenhuis realization lemma}
		Let $X$ be a $2n$-dimensional manifold and $J_x$ be an almost complex structure on the tangent space $T_xX$ at a fixed point $x\in X$.  Given a formal Nijenhuis tensor $N\colon A^{1,0}_x\to A^{0,2}_x$  there exists an almost complex structure $J$ on a neighborhood of $x$ such that $J|_{T_xX}=J_x$ with Nijenhuis tensor $N$ at $x$.
	\end{lem}
	But $\NN$ is also an affine map with contractible fibers, so it induces a bundle homotopy equivalence
	\[
	\mathcal{R}^{\geq i}\to H^{\geq i}.
	\]
	In particular, the spaces of sections are homotopy equivalent.
	
	\subsection{Applications}\label{SecApplications}
	\subsubsection{Everywhere non-integrable almost complex structures}
	Let $X$ be a $2n$-dimensional manifold and let us now discuss the special case of Theorem~\ref{thm: h-pple} where $\rank_\C N_x\geq 1$ at all points $x\in X$. Given any almost complex structure $J$ on $X$, the condition $\rank_\C (N_J)_x\geq 1$ implies the existence of a nowhere vanishing section of the bundle $H_X:=\Hom_\C(A^{1,0}_X,A^{0,2}_X)$. Equivalently, it implies that the top Chern class, say $c$, of this bundle is trivial. Conversely, if $c=0$, there exists a nowhere vanishing section $N$ of $H$. Then by Theorem~\ref{thm: h-pple} applied to the pair $(J,N)$, we see that $J$ is homotopic to an everywhere non-integrable almost complex structure $J'$. Now the bundle $H_X$ has complex rank $n\cdot\binom{n}{2}$, which is larger than the complex dimension $n$ of $X$ for $n>2$, so in this case $c=0$ necessarily. This proves:
	\begin{cor}\label{cor: eniacs}
		Any almost complex structure on a manifold of dimension $2n>4$ is homotopic to an everywhere non-integrable almost complex structure.
	\end{cor}
	\subsubsection{Maximally non-integrable structures in high dimensions}
	Let $n\geq 5$ and define $H^{deg}_X\subseteq H_X$ to be the subspace of maps which have fiberwise non-maximal rank at every point. The codimension of $H^{deg}_X$ in $H_X$ equals $\binom{n}{2} - n +1$ (c.f. \cite[2.2.1.]{eliashbergmishachev02}). Since $n\geq 5$,
	\begin{align*}
		\dim_\C\! H_X-(\dim_\C\!X + \dim_\C\!H_X^{deg}) &= \frac{1}{2}(n^2-5n+2)>0.
	\end{align*}
	Therefore, by transversality, a generic section to $H_X$ will be disjoint from $H^{deg}_X$. In other words, we can pick a global section $N\in \Gamma(X,H)$ of maximal rank. By Theorem~\ref{thm: h-pple}, there exists a homotopy of pairs from $(J,N)$ to a holonomic pair $(J',N_{J'})$ with $N_{J'}$ of maximal rank. We have shown:
	
	\begin{cor}
		For $n>4$, any almost complex structure on a $2n$-manifold is homotopic to a maximally non-integrable one.
	\end{cor}
	\begin{rem}
		The same arguments, applied to the subspace $H^{\leq r}$ of homomorphisms of pointwise rank $\leq r$, give an alternative proof of Corollary~\ref{cor: eniacs} and show that $n=4$, any almost complex structure is homotopic to one with Nijenhuis tensor of rank $\geq 3$. In particular, we have proved Corollary~\ref{corintro: novacs everywhere}.
	\end{rem}
	
	\subsubsection{$4$-manifolds}\label{Sec4Mfd}
	Let us now discuss the existence of maximally non-integrable almost complex structures in the $4$-dimensional case. The main goals of this section are the characterization of closed four manifolds admitting a maximally non-integrable almost complex structure (Corollary~\ref{corintro: dim 4}). We also determine the homotopy type of the fiber of the space of formal maximally non-integrable almost complex structures over a $4$-manifold and show that there exist homotopy classes of almost complex structures  containing infinitely many maximally non-integrable almost complex structures which are non-homotopic as maximally non-integrable almost complex structures.
	\begin{proof}[Proof of Corollary~\ref{corintro: dim 4}]
		We will prove \eqref{Eq4DimMaNIACS1} $\Longleftrightarrow$ \eqref{Eq4DimMaNIACS2} and \eqref{Eq4DimMaNIACS2} $\Longleftrightarrow$ \eqref{Eq4DimMaNIACS3}.
		
		\noindent The implication \eqref{Eq4DimMaNIACS1} $\Rightarrow$ \eqref{Eq4DimMaNIACS2} is essentially Corollary~\ref{CorObstructions}.
		Conversely, given any almost complex structure on $X$, the condition $5\chi(X)+6\sigma(X)=0$ implies that the top Chern class $c$ of the bundle $H_X$ vanishes.
		Therefore there exists a nowhere vanishing section, i.e. a formal Nijenhuis tensor of a maximally non-integrable almost complex structure. By Theorem~\ref{thm: h-pple}, this is homotopic to a genuine Nijenhuis tensor of a maximally non-integrable almost complex structure on $X$ and the first part is proved.
		
		\noindent Regarding the implication \eqref{Eq4DimMaNIACS2} $\Rightarrow$ \eqref{Eq4DimMaNIACS3}, recall that $4\vert(\chi+ \sigma)$ holds on any almost complex $4$-manifold. Hence we have the following congruences modulo $4$
		\begin{equation}\label{EqCongruence}
			0\equiv 5\chi(X)+6\sigma(X)\equiv 5(\chi(X)+ \sigma(X))+\sigma(X)\equiv \sigma(X).
		\end{equation}
		Thus we are left needing to show that an almost complex $4$-manifold satisfying $5\chi(X)+6\sigma(X)=0$ cannot have positive definite intersection form and $b_2(X)>0$.
		An almost complex $4$-manifold with positive definite intersection form satisfies $\chi(X)+\sigma(X)\geq0$ and $4\vert(\chi(X)+ \sigma(X))$, see \cite[Theorem~8.B]{matsushita91}. Then we have
		$$-\sigma(X)=5\left(\chi(X)+\sigma(X)\right)\geq0.$$
		Thus, a $4$-manifold $X$ satisfying \eqref{Eq4DimMaNIACS2} can have positive definite intersection form only if $\sigma(X)=b_2(X)=0$. 
		
		\noindent Conversely, assuming \eqref{Eq4DimMaNIACS3} we want to show that $X$ admits an almost complex structure. We will now distinguish three cases:
		\begin{enumerate}[(i)]
			\item $b_2(X)=0$. \label{EqCaseB20}
			\item The intersection form of $X$ is indefinite.\label{EqCaseIndefinite}
			\item The intersection form of $X$ is negative definite.\label{EqCaseNegDef}
		\end{enumerate}
		If $b_2(X)=0$, then $X$ admits an almost complex structure if $2\chi(X)+3\sigma(X)=0$, cf. for instance \cite[Theorem~8.B]{matsushita91}. But this is immediate since $b_2(X)=0$ implies $\sigma(X)=0$ and this, together with $0=5\chi(X)+6\sigma(X)$, yields $\chi(X)=0$. This settles the case \eqref{EqCaseB20}.
		For \eqref{EqCaseIndefinite}, combining $5\chi(X)+6\sigma(X)=0$ with $4\vert\sigma(X)$ we get $\chi(X)+\sigma(X)\equiv0 \pmod{4}$. This last congruence is known to be a sufficient condition for a closed $4$-manifold with indefinite intersection form to admit an almost complex structure, see for example \cite[Theorem~8.A]{matsushita91}.
		We are left with the case \eqref{EqCaseNegDef}. Since $4\vert\sigma(X)$ (which translates to $4\vert b_2(X)$), we can restrict to the case $b_2(X)\geq4$.
		In this instance $X$ admits an almost complex structure if the congruence  $\chi(X)+\sigma(X)\equiv0 \pmod{4}$ is satisfied \cite[Theorem~8.C1]{matsushita91}.
		This congruence, as in the previous case, follows by combining $5\chi(X)+6\sigma(X)=0$ with $4\vert\sigma(X)$.
	\end{proof}
	
	Let us now discuss the homotopy type of the fiber of the space $H_X^{mni}$ of formal maximally non-integrable almost complex structures on a $4$-manifold $X$.
	Let $\pi\colon Z\lra X$ be the twistor space of $X$.
	Recall from Section~\ref{NijenhuisTensors} that the space of linear Nijenhuis tensors $H_X$ is defined as $A^{0,2}\otimes T^{1,0}$, see \eqref{EqNijenhuisBundle}. In order to study its homotopy type we must first understand the homotopy type of the normal bundle to the fibers $F$ of $Z$, i.e. the restriction of $T$ to $F$. Since $F$ is homotopy equivalent to $S^2$, $T_{\vert F}$ is homotopy equivalent to a rank $2$ complex vector bundle on $S^2$. This bundle is $L\oplus L$ where $L$ has Euler class $x$ the positive generator of $H^2(S^2)$, see  \cite[page~437]{AHS78}. Therefore $c_1(H_X|_F)\cong 6x$. Since $H_X^{mni}$ is just $H_X$ without the zero section, its fiber has the homotopy type of the sphere bundle associated with $H|_{F}$, which is an $S^3$-bundle over $S^2$. Since $2\mid c_1(H|_F)$, this sphere bundle is homotopy equivalent to the trivial bundle $S^2\times S^3$.
	
	We now specialize this discussion to the case of parallelizable $4$-manifolds. In particular we obtain the following
	
	\begin{prop}\label{PropPar4Mfd}
		Let $X$ be a parallelizable $4$-manifold equipped with an arbitrary almost complex structure $J$. Then homotopy classes of maximally non-integrable almost complex structures which are homotopic to $J$ as almost complex structures are in one-to-one correspondence to homotopy classes $[X,S^3]$ of maps $X\lra S^3$.
	\end{prop}
	\begin{proof}
		Let $J$ be an arbitrary almost complex structure on a parallelizable $4$-manifold $X$. 
		This corresponds to a section $J\colon X\lra Z$ of the twistor bundle $Z=X\times F$. 
		Recall that $Z$ is homotopy equivalent to the trivial $S^2\sqcup S^2$-bundle over $X$ by triviality of $TX$.
		Homotopy classes of maximally non-integrable almost complex structures in the same homotopy class of $J$ are the same as homotopy classes of (nowhere-vanishing) extensions of $J$ to $H$, i.e. (nowhere-vanishing) sections of the bundle $J^*H$.
		We claim that the bundle $J^*H$ over $X$ is trivial. 
		In order to see this recall that a rank $4$ real vector bundle over a closed $4$-manifold is determined by its second Stiefel-Whitney class $w_2$, its first Pontryagin class $p_1$ and its Euler class $e$. 
		Since $J^*H$ has a complex structure induced by $J$, we have $c_1\equiv w_2 \pmod 2$, $p_1=c_1^2-2c_2$ and $e=c_2$. 
		Now, using the fact that $J^*H\cong\Hom\left(A_X^{0,2}\otimes TX^{1,0}\right)$, we get 
		\begin{enumerate}[$\bullet$]
			\item $w_2(J^*H)\equiv c_1(J^*H)=3c_1(X)\equiv c_1(X)\equiv 0 \pmod 2$.
			\item $e(J^*H)=5\chi+6\sigma=0$ where the first equality follows from Corollary~\ref{CorObstructions} and the second one follows from $\chi=\sigma=0$.
			\item $p_1(J^*H)=c_1^2(J^*H)-2c_2(J^*H)=0$ by the previous points and $c_1^2(X)=0$.
		\end{enumerate}
		Therefore maximally non-integrable almost complex structures in the same homotopy class of $J$ are maps $X\lra S^3$ because $J^*H\cong X\times S^3$. Consequently homotopy classes of maximally non-integrable almost complex structures in the same homotopy class of $J$ correspond to elements of $[X,S^3]$. This concludes the proof.
	\end{proof}
	
	By picking $X$ so that there are infinitely many homotopy classes of maps $X\lra S^3$ we get:
	\begin{cor}
		There exist parallelizable $4$-manifolds for which each homotopy class of almost complex structures contains infinitely many maximally non-integrable almost complex structures which are non-homotopic as maximally non-integrable almost complex structures.
	\end{cor}
	\begin{proof}
		Let $X$ be the parallelizable $4$-manifold $S^1\times S^3$ equipped with an arbitrary almost complex structure $J$. Consider the maps $f_n\colon X\lra S^3$ which are constant in the $S^1$ factor and have degree $n$ in the  $S^3$ factor $ S^3$.
		By Proposition~\ref{PropPar4Mfd} these maps represent maximally non-integrable almost complex structures on $X$ which are not homotopic as maximally non-integrable almost complex structures even though they are all homotopic to $J$ as almost complex structures.
	\end{proof}
	
	\subsection{6-manifolds}
	Let $X$ be an oriented 6-manifold. Assume that $X$ admits a maximally non-integrable almost complex structure, then by Corollary \ref{CorObstructions} its first Chern class $c_1$ is 3-torsion, i.e. $3c_1=0$. Moreover, $c_1$ is an integral lift of the second Stiefel-Whitney class $w_2$, which is 2-torsion. Therefore $w_2=0$, that is, $X$ is $\mathrm{spin}$.
	
	Conversely, if $w_2=0$, then we may take $0$ as an integral lift and there is an almost complex structure $J$ (unique up to homotopy) such that $c_1=0$ (cf. \cite[Thm 9]{wall66}). The condition $c_1=0$ is equivalent to $A_X^{3,0}$ being the trivial complex line bundle, i.e. there exists a nowhere vanishing $(3,0)$-form $\eta$. Picking a Hermitian metric yields a $\C$-linear Hodge operator $\star:A^{1,0}_X\to A^{3,2}_X$. Combining these two we obtain a bundle isomorphism $N:=(\cdot\wedge\eta)^{-1}\circ\star:A^{1,0}_X\to A^{0,2}_X$. By {Theorem \ref{thm: h-pple}} applied to the pair $(J,N)$, the almost complex structure $J$ is homotopic to a maximally non-integrable one.
	
	Note that one can always kill torsion classes in $H^2(X,\Z)$ by passing to a finite cover. Therefore, if for a given almost complex structure $J$ on a $6$-manifold $X$ we have that $c_1$ is torsion, then there exists an almost complex structure on a finite cover of $X$ with $c_1=0$. We have thus proved Corollary \ref{corintro: dim 6}

	\subsubsection{Parallelizable manifolds}
	On a $2n$-dimensional manifold $X$ with trivial tangent bundle $TX$ the choice of a frame $V_1,...,V_{2n}$ induces a trivialization of the bundle $Z\cong X\times F\to X$, where $F$ is any fiber of $\pi:Z\to X$. One may define an almost complex structure via $JV_{2i-1}=V_{2i}$ for all $i=1,...,n$. This corresponds to a constant section $X\to X\times F$ gives a complex trivialization of the tangent bundle. Conversely, for every almost complex structure which makes $TX$ into a trivial complex vector bundle, we may find a frame in which it becomes constant as above. In such a situation, the bundle $H_X$ is trivial as a complex vector bundle and so one may pick a section of maximal rank. Hence, we have proved Corollary~\ref{corintro: parallelizable}
	
	\subsection{Examples}
	\subsubsection{Mapping tori}
	This example follows the idea presented by Geiges in \cite{geiges00}.
	Let $X$ be a parallelizable $4$-manifold which fibers over $S^1$, that is $\pi\colon X\lra S^1$ is a bundle with fiber a (oriented, closed) $3$-manifold $Y$.
	In other terms, $X$ is the suspension of an orientation preserving diffeomorphism $f\colon Y\lra Y$.
	In particular, $X$ admits a maximally non-integrable almost complex structure by Theorem~\ref{corintro: parallelizable}. We construct now an explicit one. \\
	
	Denote by $T$ the nowhere vanishing vector field on $X$ defined by $\pi_*T=\partial_t$ where $t$ is the coordinate on $S^1$. Then we can always complete $T$ to a trivialization $\{T,U,V,Z\}$ of $TX$, for instance taking $U=iT$, $V=jT$ and $Z=kT$ under the identification $TX\cong X\times\HH$.
	From this we construct the following trivialization of $TX$.
	\begin{align*}
		A_n&=T+\frac{1}{n}\sin\left(n^2\pi^*t\right)U-\frac{1}{n}\cos\left(n^2\pi^*t\right)V\\
		B_n&=Z+\frac{1}{n}\cos\left(n^2\pi^*t\right)U+\frac{1}{n}\sin\left(n^2\pi^*t\right)V\\
		C_n&=\frac{1}{n}U\\
		D_n&=\frac{1}{n}V.
	\end{align*}
	Now define an almost complex structure $J$ on $X$ by setting
	$$JA_n=C_n\hspace{3mm}\mbox{and} \hspace{3mm}JB_n=D_n\,\,.$$
	Then $J$ is maximally non-integrable (for $n$ large enough) because the Nijenhuis tensor $N_J$ is everywhere non-trivial.
	In particular the vector field $N_J\left(A_n,B_n\right)$ is nowhere vanishing. 
	In order to show this we begin by noticing that 
	\begin{align*}
		\frac{1}{n}N_J\left(A_n,B_n\right)=&\frac{1}{n^3}[U,V]-\frac{1}{n}[A_n,B_n]-\frac{1}{n^2}J[U,B_n]-\frac{1}{n^2}J[A_n,V]\\
		=&-\frac{1}{n}[A_n,B_n]+O(n^{-1}).
	\end{align*}
	Using $\lie_{T}\pi^*t=1$ and setting $a\colon=\lie_Z\pi^*t$ we have
	\begin{align*}
		-\frac{1}{n}[A_n,B_n]=&\left(\sin(n^2\pi^*t)+a\cos(n^2\pi^*t)\right)U\\
		&-\left(\cos(n^2\pi^*t)-a\sin(n^2\pi^*t)\right)V+O(n^{-1})
	\end{align*}
	which is nowhere-vanishing for $n$ large enough.
	We conclude that $N_J\neq 0$ for all $x\in X$ so that $J$ is maximally non-integrable. 
	
	\subsubsection{Homogeneous manifolds}
	Consider homogeneous manifolds $X=G/\Gamma$ for a Lie group $G$ and a discrete subgroup $\Gamma$, equipped with a left-invariant complex structure $J$. By their homogeneous nature, the Nijenhuis tensor $\mubar_J$ always has constant rank. Thus, it suffices to check that $\mubar_J$ has maximal rank at a point. For instance, in real dimension $4$, any non-integrable left-invariant complex structure is maximally non-integrable. Further, one has a complex trivialization 
	\begin{equation}\label{eqn: H-triv}
		\Hom(A^{1,0}_X,A^{0,2}_X)\cong\underline{\C}\otimes \Hom(\fg^{1,0},\fg^{0,2}),
	\end{equation} where $\fg^{p,q}$ denotes the space of left-invariant $(p,q)$-forms. Therefore, by Theorem~\ref{thm: h-pple}, resp. Corollary~\ref{corintro: parallelizable}, any left-invariant almost complex structure is homotopic to a maximally non-integrable one. However, the latter need not be left-invariant anymore, as the following (well-known) observation shows:
	
	\begin{obs}
		Every left-invariant almost complex structure on the torus $T^{2n}$ is integrable.
	\end{obs}
	
	In fact, this is an immediate consequence of the trivialization \eqref{eqn: H-triv} and the fact that $\de$ (hence a fortiori $\mubar$) vanishes on left-invariant forms.
	
	\begin{remark}
		In dimension $4$, $T^4$ is in fact the only solvmanifold which does not admit a left-invariant maximally non-integrable almost complex structure. 
	\end{remark}
	
	Let us check constructively, and without using Theorem~\ref{thm: h-pple} that there exist maximally non-integrable almost complex structures on all tori. The idea becomes apparent in dimension $4$ and for simplicity we consider the complex plane $X=\C^2$ first. 
	Let $\omega_1:=\de z_1+f\de \bar z_2,\omega_2:=\de z_2$ for some function $f$ and define an almost complex structure by setting $A^{1,0}:=\langle\omega_1,\omega_2\rangle$. Then $A_X^{0,2}=\langle\bar\omega_1\wedge\bar\omega_2\rangle$ and hence \begin{align*}
		\mubar(\omega_1)
		&=\pr^{0,2}\bigg(\frac{\del}{\del z_1}f \de z_1\wedge \de\bar z_2 + \frac{\del}{\del z_2} f \de z_2\wedge \de\bar z_2 + \frac{\del}{\del \bar z_1}f \de\bar z_1\wedge \de \bar z_2\bigg)\\
		&=\pr^{0,2}\bigg(\frac{\del}{\del z_1}f (\omega_1-f\bar\omega_2)\wedge \bar\omega_2 + \frac{\del}{\del z_2} f \omega_2\wedge \bar \omega_2 + \frac{\del}{\del \bar z_1}f (\bar\omega_1-\bar f\omega_2)\wedge  \bar \omega_2\bigg)\\
		&=\frac{\del}{\del \bar z_1} f\bar\omega_1\wedge\bar\omega_2.
	\end{align*}
	Hence, for any $f$ such that $\frac{\del}{\del \bar z_1}f$ is nowhere vanishing, the almost complex structure is maximally non-integrable. 
	For a $4$-torus $\C^2/(2\pi\Z[i]^2)$, one can pick for instance. $f= \exp(iRe(z_1))$.
	
	Generalizing these ideas, we obtain:
	\begin{thm}
		On every even dimensional torus $T^{2n}$, $n\geq 2$, there exists a maximally non-integrable almost complex structure.
	\end{thm}
	\begin{proof}
		The case $n=2$ being treated above, it is enough to restrict to the case $n\geq 3$. Let $T:=\C^n/(2\pi \Z[i]^n)$ with complex coordinates $z_1,...,z_n$. Define an almost complex structure by deforming the space of $(1,0)$-forms as follows.
		Consider a basis of $(1,0)$-forms given by
		\[
		\omega_i:= \de z_i+a_i \de\bar z_i
		\]
		where
		\[
		a_i:= f(z_{i+1}) + A
		\]
		with $f:\C/(2\pi\Z[i])\rightarrow S^1$ given by $f(z):=\exp(iRe(z))$, $A\in\C$ a constant s.t. $\|A\|\geq 2$ and by abuse of notation $z_{n+1}:=z_1$. Note that we have
		\begin{align*}
			\frac{\del}{\del z_i}a_j&= 0&i\neq j+1\\
			\frac{\del}{\del z_i}a_j&= \frac{i}{2}f(z_{i+1})&i=j+1\\
			\frac{\del}{\del \bar z_i}a_j&= 0&i\neq j+1\\
			\frac{\del}{\del \bar z_i}a_j&= -\frac{i}{2}f(z_{i+1})&i=j+1
		\end{align*}
		where again indices are to be read modulo $n$. We can re-express the $\de z_i$ from the $\omega_i$ as
		\[
		\de z_i=b_i(\omega_i-a_i\bar\omega_i)
		\]
		with $b_i:=(1-\|a_i\|^2)^{-1}$. By assumption on $A$, we have $\|a_i(z)\| >1$ at every point $z\in T$ so that $b_i$ never vanishes.
		\begin{align*}
			\mubar(\omega_i)
			&=\pr^{0,2}\bigg(\frac{\del}{\del z_{i+1}}a_i \de z_{i+1}\wedge \de \bar z_i + \frac{\del}{\del\bar z_{i+1}}a_i \de \bar z_{i+1}\wedge \de\bar z_i\bigg)\\
			&=\underbrace{\frac{i}{2}b_ib_{i+1}f(z_{i+1})(1-a_{i+1})}_{\neq 0}\bar\omega_{i+1}\wedge\bar\omega_i.
		\end{align*}
		Since a basis for $A^{0,2}_{T}$ is given by the $\bar\omega_i\wedge\bar\omega_j$ with $i>j$, all $\mubar(\omega_i)$ are linearly independent, i.e. $\mubar$ has maximal rank.
	\end{proof}
	\begin{remark}
		Taking products one deduces that for any $n,d$ with $0\leq d\leq \min\{n,\frac{n(n-1)}2\}$, there exists an almost complex structure on $T^{2n}$ with Nijenhuis tensor of constant rank $d$.\footnote{Strictly speaking, products of the previous examples cover all cases except a rank $2$ structure on $T^6$. Such a structure can be obtained for example by defining $\omega_1:=\de z_1 + \exp(i Re(z_1))\de\bar{z}_3$, $\omega_2:=\de z_2 + \exp(i Re(z_2))\de\bar{z}_3$ and $\omega_3:=\de z_3$.}
	\end{remark}
	
	\section{Dolbeault cohomology}\label{SecDolbeault}
	\subsection{Preliminaries} We recall the notion of Dolbeault cohomology introduced in \cite{ciriciwilson18}: On any almost complex manifold $(X,J)$, we write $d=\mubar+\delbar+\del+\mu$ where the summands are the components of bidegrees $(-1,2)$, $(0,1)$, $(1,0)$ and $(2,-1)$. In particular, the $\mubar$ from the previous section is the restriction of the $\mubar$ from this section to $A^{1,0}_X$. The equation $\de^2=0$ translates into:
	\begin{align*}
		0&=\mubar^2\\
		0&=\mubar\delbar+\delbar\mubar\\
		0&=\mubar\del+\del\mubar+\delbar^2\\
		0&=\mubar\mu+\delbar\del+\del\delbar+\mu\mubar\addtocounter{equation}{1}\tag{\theequation}\label{4-complex}\\
		0&=\mu\delbar+\delbar\mu+\del^2\\
		0&=\mu\del+\del\mu\\
		0&=\mu^2.
	\end{align*}
	By the first three of these equations, the following generalization of Dolbeault-cohomology is well-defined:
	\[
	H_{Dol}(X):=H_{Dol}(X,J):=H_{\delbar}(H_\mubar(\cA_X^{\bullet,\bullet})).
	\]
	\subsection{Proof of Theorem~\ref{thm: mani => HDol infty}}
	We prove first a technical lemma. Denote by $\cA^1_{X,\R}:=\Gamma(X,TX^\vee)$ the space of (real valued) $1$-forms.
	\begin{lem}\label{lem: quot inf}
		Let $X$ be any smooth manifold of dimension $n\geq 3$ and $E\in \Gamma(X,\operatorname{End}(TX))$ an endomorphism of the tangent bundle. Then the $\R$-vector space
		\[(\im d + \im E\circ d)\cap \cA^1_{X,\R}
		\]
		has infinite codimension in $\cA_{X,\R}^1$. More precisely, at any given point $x\in X$, there exists a $k\geq 1$ such that the space $(\im d + \im E\circ d)^{(k)}$ of $k$-jets at $x$ of elements in $(\im d + \im E\circ d)^{(k)}$ has positive codimension in the space $(\cA_{X,\R}^1)^{(k)}$ of all $k$-jets of $1$-forms at $x$.
	\end{lem}
	\begin{proof}
		We can assume $X$ to be a disk around $x=0$ with coordinates $x_1,...,x_n$. The space $(\cA_{X,\R}^0)^{(k)}$ of $k$-jets of functions ($k$-th order Taylor expansions) at $x$ may be identified with the space of polynomials of degree $\leq k$ in the $x_i$, therefore
		\[
		\dim (\cA^0_{X,\R})^{(k)}= \binom{n+k}{k}
		\]
		Hence,
		\[
		\dim (\cA_{X,\R}^1)^{(k)}= n\cdot\binom{n+k}{k}
		\]
		and 
		\[
		\dim (\im d + \im E\circ d)^{(k)}\geq 2\cdot \binom{n+k+1}{k+1}.
		\]
		But 
		\[
		n\cdot\binom{n+k}{k}> 2\cdot \binom{n+k+1}{k+1}\Leftrightarrow 0> n+(2-n)k+2,
		\]
		and since $n\geq 3 $ there is always a $k$ satisfying this.
	\end{proof}
	
	\begin{proof}[Proof of Theorem~\ref{thm: mani => HDol infty}]
		Let us first assume that $J$ is everywhere maximally non-integrable. Then the Dolbeault cohomology in degree $(0,1)$ is the cokernel $H_{Dol}^{0,1}(X)=\cA^{0,1}_X/(\im\delbar)^{0,1}$.  There is an $\R$-linear isomorphism
		\begin{align*}
			\cA^{0,1}_X/(\im\delbar)^{0,1}\rightarrow& \cA^1_{X,\R}/(\im d + J\im d)\\
			\omega\mapsto&\omega + \bar\omega. 
		\end{align*}
		and the right hand side is infinite dimensional by Lemma~\ref{lem: quot inf}. 
		
		For the general case, let $x\in X$ be a point s.t. $J_x$ is maximally non-integrable. Pick an open neighborhood $x\in U\subseteq X$ such that $J|_U$ is maximally non-integrable, so that $H_{Dol}^{0,1}(U)=\cA^{0,1}_{U}/\im \delbar$. By Lemma~\ref{lem: quot inf}, we may pick an infinite dimensional subspace $K\subseteq \cA^{0,1}_U$ consisting of forms with compact support s.t. $K\rightarrow H_{Dol}^{0,1}(U)$ is injective. Extending these forms by zero and projecting yields a well-defined map 
		\[
		K\rightarrow H_{Dol}^{0,1}(X)=\frac{\delbar^{-1}\im \mubar\cap \cA_X^{0,1}}{\im\delbar\cap \cA_X^{0,1}}.
		\]
		The result now follows from the commutativity of the following diagram:
		\[
		\begin{tikzcd}
			K\ar[r]\ar[rd, hook]&H_{Dol}^{0,1}(X)\ar[d]\\
			&H_{Dol}^{0,1}(U).
		\end{tikzcd}
		\]
	\end{proof}
	
	\begin{rem}
		A spectral sequence argument shows that in the setting of Theorem \ref{thm: mani => HDol infty}, also $H_{Dol}^{1,1}(X)$ is infinite dimensional.
	\end{rem}
	\begin{ex}
		Let $T$ be the 4-torus endowed with the maximally non-integrable almost complex structure constructed above. Then the Dolbeault cohomology $H_{Dol}^{0,1}(T)$ is infinite dimensional as a $\C$-vector space.
	\end{ex}
	\begin{proof}
		Since the almost complex structure is maximally non-integrable, $H_{Dol}^{0,1}(T)=\cA^{0,1}_T/\im \delbar$. Now, for any function $g\in \cA^{0,0}_T$, we have
		\begin{align*}
			\delbar g&=\pr^{0,1}\bigg(\frac{\del}{\del z_1} g\de z_1+ \frac{\del}{\del z_2}g \de z_2 +\frac{\del}{\del \bar z_1} g\de \bar z_1+ \frac{\del}{\del \bar z_2}g \de\bar z_2\bigg)\\
			&=\pr^{0,1}\bigg(\frac{\del}{\del z_1}g(\omega_1-f\bar\omega_2)-\frac{\del}{\del z_2}g\omega_2+\frac{\del}{\del \bar z_1} g (\bar\omega_1-\bar f\omega_2)+ \frac{\del}{\del \bar z_2}g \bar\omega_2\bigg)\\
			&=\frac{\del}{\del \bar z_1}g\bar\omega_1+(\frac{\del}{\del \bar z_2}g-f\frac{\del}{\del z_1}g)\bar\omega_2
		\end{align*}
		That is, if an element of the form $G\bar\omega_2$ with $G\in \cA^{0,0}_T$ represents the zero-class in $H^{0,1}_{Dol}(T)$, there exists $g\in \cA^{0,0}_T$ such that
		\begin{align}
			\frac{\del}{\del \bar z_1}g&=0\label{Eq1Torus}\\
			\frac{\del}{\del \bar z_2}g-f\frac{\del}{\del z_1}g&=G\, .\label{Eq2Torus}
		\end{align}
		The first equation says that $g$ is holomorphic in $z_1$. Since $T$ is compact, $g$ is constant on each slice of the type $(z_1,x) $ for fixed $x$, i.e. $g$ is independent of $z_1$. Therefore, also $\frac{\del}{\del z_1} g=0$ and \eqref{Eq2Torus} simplifies to
		\[
		\frac{\del}{\del \bar z_2}g=G\, .
		\]
		This, together with \eqref{Eq1Torus}, implies $\frac{\del}{\del\bar z_1}G=0$, hence also $G$ has to be independent of $z_1$. 
		
		Let $B$ be a linear subspace of $\cA^{0,0}$ which is contained in the complement of the space of functions which are independent of $z_1$. For such $B$ the subspace $\cA^{0,1}_B:=\{ G\bar\omega_2\mid G\in B\}$ injects into $H_{Dol}^{0,1}(T)$ by our previous argument. Now choose $B$, and hence $\cA_B^{0,1}$, to be infinite dimensional.
	\end{proof}
	\begin{rem}\label{rem: torus n,n-1}
		A similar calculation shows that $H_{Dol}^{2,1}(T)$ is infinite dimensional for this almost complex structure.
	\end{rem}
	\begin{rem}[Dolbeault cohomology is not Hausdorff]
		Consider the sequence of functions 
		\[
		\tilde{g}_n(z_1,\bar{z}_1,z_2,\bar{z_2}):=\frac{\bar{z}_1}{n} + z_1\in \cA^0_{\C^2}.
		\]
		Let $U\subseteq \C^2$ be an open set for which $\pr: U \to T^4$ is an embedding. Pick $V\subseteq U$ and a smooth function $\rho:U\to\R$ s.t. $\rho|_V\equiv 1$ and $\rho|_{\C^2\setminus U}\equiv 0$. Then we may consider $g_n:=\rho\cdot \tilde{g}_n$ as functions on the torus. We claim that with respect to the previously considered almost complex structure on $T^4$, the limit of the $\delbar g_n$ for $n\to\infty$ does not lie in the image of $\delbar$. Indeed, on $V$ we have
		\[
		\delbar g_n\overset{n\to\infty}\longrightarrow f\omega_2\not\in\im\delbar.
		\]
		In other words, the zero class in $H_{Dol}^{0,1}(T^4)$ is not closed, so $H_{Dol}^{0,1}(T^4)$ is not Hausdorff. A natural question, suggested to us by Jean Ruppenthal, is whether for compact $X$, the maximal Hausdorff quotient $(H_{Dol}^{0,1}(X))^{HD}:=H_{Dol}^{0,1}(X) / \overline{[0]}$ is finite dimensional. Note however that for maximally non-integrable $4$ or $6$-manifolds, we have $H_{Dol}^{n,n-1}(X)=\ker\delbar$, which is Hausdorff, but can still be infinite dimensional as Remark \ref{rem: torus n,n-1} shows.
	\end{rem}
	\appendix
	\section{Bott-Chern and Aeppli cohomologies}\label{AppendixBottChern}
	We recall an algebraic notion which can be seen as an abstraction of the complex of forms on an almost complex manifold:
	\begin{definition}
		A complex $(A,d)$ where $A$ carries a bigrading and $d$ has bidegree $|d|=(-1,2)+(0,1)+(1,0)+(2,-1)$ is called a $4$-complex.
	\end{definition}
	In the following, $A$ will always denote a $4$-complex and we will denote the components of $d$ by $\mubar,\delbar,\del,\mu$ with $|\mubar|=(-1,2)$, $|\delbar|=(0,1)$, $|\del|=(1,0)$ and $|\mu|=(2,-1)$. As before, the property $d^2=0$ translates into the equations \eqref{4-complex}. For technical reasons, we will assume all our $4$-complexes to be bounded, i.e. $A^{p,q}=0$ for all but finitely many $(p,q)$. Following \cite{ciriciwilson18}, for any such $4$-complex one may introduce the Dolbeault and `conjugate' Dolbeault cohomology of $A$:
	\[
	H_{Dol}(A):=H_{\delbar}(H_{\mubar}(A))\qquad H_{Dolbar}(A):=H_{\del}(H_{\mu}(A)).
	\]
	(which is meaningful by the equations \eqref{4-complex}) and there will be the following spectral sequences, functorial for maps of $4$-complexes: (see \cite{ciriciwilson18})
	\[
	H_{Dol}(A)\Longrightarrow H_d(A)\Longleftarrow H_{Dolbar}(A).
	\]
	A $4$-complex with $\mubar=\mu=0$ is called a double complex. If $A$ is a double complex one may define the Bott-Chern and Aeppli cohomologies:
	
	\[
	H_{BC}(A):=\frac{\ker\del\cap\ker\delbar}{\im\del\delbar}\qquad H_{A}(A):=\frac{\ker\del\delbar}{\im\del+\im\delbar}.
	\]
	These fit into a commutative diagram
	\[
	\begin{tikzcd}\label{cohomology diagram}\tag{$\ast$}
		&H_{BC}(A)\ar[ld]\ar[d]\ar[rd]&\\
		H_{\delbar}(A)\ar[rd]\arrow[Rightarrow]{r}&H_{d}(A)\ar[d]&H_{\del}(A)\ar[ld]\arrow[Rightarrow]{l}\\
		&H_A(A).
	\end{tikzcd}
	\]
	Note that for a double complex $H_{Dol}(A)=H_\delbar(A)$, so the middle row is meaningful for any $4$-complex. Our goal will be to extend the definition of Bott-Chern and Aeppli cohomology to almost complex manifolds or, more generally, arbitrary $4$-complexes
	\begin{definition}
		Given a $4$-complex $(A,\mubar,\delbar,\del,\mu)$, we say
		\[
		A_s:=\ker\mubar\cap\ker\delbar^2\cap\ker\del^2\cap\ker\mu
		\]
		is the sub double complex and
		\[
		A_q:= A/(\im\mubar+\im\delbar^2+\im\del^2+\im \mu)
		\]
		the quotient double complex.
	\end{definition}
	
	\begin{lem}
		The maps $\del$ and $\delbar$ induce well-defined maps on $A_s$ and $A_q$. With these, $(A_s,\del,\delbar)$ and $(A_q,\del,\delbar)$ are double complexes.
	\end{lem}
	\begin{proof}
		We only do $A_s$, the other case is essentially the same. First of all, observe that $A_s$ is still bigraded, since it is the intersection of kernels of maps that have a single bidegree. The equations of a $4$-complex simplify on the subspace $A_s$ to
		\begin{equation}\label{eqn: A_s eqns}
			0=\mubar\delbar=\mubar\del=\delbar\del+\del\delbar=\mu\delbar=\mu\del
		\end{equation}
		
		Given $\omega\in A_s$ let us check that $\del\omega\in A_s$: We have $\mubar\del\omega=\mu\del\omega=0$ by the equations \eqref{eqn: A_s eqns} and $\del^2\del\omega=\del\del^2\omega=0$ since $A_s\subseteq \ker\del^2$. Finally $\delbar^2\del\omega=-\mubar\del^2\omega-\del\mubar\del\omega=0$. Analogously, $\delbar\omega\in A_s$. Finally the equation $\del\delbar+\delbar\del=0$ from \eqref{eqn: A_s eqns}, together with $\del^2=\delbar^2=0$ show that we are in a double complex.
	\end{proof}
	
	Since $A_s$ and $A_q$ are double complexes, the usual definition of Bott-Chern and Aeppli cohomology are meaningful for them. So we get diagrams
	\[\begin{tikzcd}
		&H_{BC}(A_s)\ar[ld]\ar[d]\ar[rd]&\\
		H_{\delbar}(A_s)\arrow[Rightarrow]{r}&(H_{dR}(A_s),F,\bar F)&H_{\del}(A_s)\arrow[Rightarrow]{l}
	\end{tikzcd}\]
	and
	\[\begin{tikzcd}
		H_{\delbar}(A_q)\ar[rd]\arrow[Rightarrow]{r}&(H_{dR}(A_q),F,\bar F)\ar[d]&H_{\del}(A_q)\ar[ld]\arrow[Rightarrow]{l}\\
		&H_A(A_q).
	\end{tikzcd}\]
	
	If we started with a double complex $A$, we have $A_s=A=A_q$ and we re-obtain the diagram $(\ast)$. It remains to connect these diagrams with the actual $A$ for a proper $4$-complex:
	
	By construction one has two maps of 4-complexes (where the outer ones have trival $\mu$ and $\mubar$) 
	\[
	A_s\longrightarrow A\longrightarrow A_q.
	\]
	Applying the de-Rham cohomology-functor, we get maps
	\[
	H_{dR}(A_s)\longrightarrow H_{dR}(A)\longrightarrow H_{dR}(A_q).
	\]
	Applying $H_{Dol}$ and using that $\mubar=0$ on $A_s$ and $A_q$, we obtain
	\[
	H_\delbar(A_s)=H_{Dol}(A_s)\longrightarrow H_{Dol}(A)\longrightarrow H_{Dol}(A_q)=H_{\delbar}(A_q)
	\]
	In summary, we have:
	\begin{prop}\label{prop: cohomology diagram}
		Given a $4$-complex $A$, set $H_{BC}(A)=H_{BC}(A_s)$ and $H_{A}(A):=H_A(A_q)$. Then there is a commutative diagram:
		\[
		\begin{tikzcd}
			&H_{BC}(A)\ar[ld]\ar[d]\ar[rd]&\\
			H_{Dol}(A)\ar[rd]\arrow[Rightarrow]{r}&(H_{dR}(A),F,\bar F)\ar[d]&H_{Dolbar}(A)\ar[ld]\arrow[Rightarrow]{l}\\
			&H_A(A)
		\end{tikzcd}
		\]
	\end{prop}
	Concerning multiplicative structures, we have proved:
	\begin{prop}
		If $A$ carries a product $\wedge$ making $(A,\wedge, d)$ into a bigraded differential algebra, $A_s$ is subalgebra and $A_q$ is a bigraded (bi)differential module over $A_s$. Hence, $H_A(A)$ is a bigraded module over $H_{BC}(A)$. 
	\end{prop}
	\begin{proof}
		Note that for $\delta\in\{\mubar,\delbar^2,\del^2,\mu\}$ and pure-type forms $\omega,\omega'\in A$, we have 
		\[
		\delta(\omega\wedge\omega')=\delta\omega\wedge\omega' \pm \omega\wedge\delta\omega'.
		\]
		From this, it is immediate that $A_s$ is closed under $\wedge$ and that $A_s \wedge (\im\mubar+\im\delbar^2+\im\del^2+\im\mu)\subseteq (\im\mubar+\im\delbar^2+\im\del^2+\im\mu)$, i.e. $A_q$ is a module over $A_s$.
		The final statement follows from the more general fact that for any bigraded differential module $(M,\del,\delbar)$ over a bigraded bidifferential algebra $(B,\del,\delbar)$, the Aeppli cohomology $H_A(M)$ is a bigraded module over the Bott-Chern cohomology $H_{BC}(B)$. In fact, let $b\in \ker \del\cap\ker\delbar \cap B^{p,q}$ and $m\in \ker\del\delbar\cap M^{r,s}$. Then 
		\[
		\del\delbar (b\wedge m)=\del\delbar b\pm \delbar b\wedge\del m \mp \del b\wedge\delbar m + b\wedge\del\delbar m=0.
		\]If additionally $b=\del\delbar b'$, then 
		\[b\wedge m=\del\delbar(b'\wedge m)\pm \del(b'\wedge \delbar m)\mp \delbar(b'\wedge\del m)\in (\im\del+\im\delbar).\] Similarly, if $m=\del m'+\delbar m''$, then $b\wedge m=\pm \del(b\wedge m')\pm \del(b\wedge m'')\in (\im\del+\im\delbar)$.
	\end{proof}
	\begin{rem}
		The statement about the module structure implies in particular that for all $p,q,r,s\in\Z$, the wedge-product induces well-defined `duality' pairings
		\[
		H_{BC}^{p,q}(A)\times H_A^{r,s}(A)\longrightarrow H_A^{p+r,q+s}(A).
		\]
	\end{rem}
	For completeness, we make the following definition, which the reader will have guessed by now:
	\begin{definition}
		Let $X$ be an almost complex manifold. The Bott-Chern and Aeppli cohomologies are defined to be
		\[
		H_{BC}(X):=H_{BC}(\cA_X)\quad H_{A}(X):=H_A(\cA_X).
		\]
	\end{definition}
	With this definition we obtain:
	\begin{thm}
		Bott-Chern and Aeppli cohomology give contravariant functors
		\begin{align*}
			H_{BC}:&\{\text{almost complex manifolds}\}\to\{\text{bigraded-commutative }\C\text{-algebras}\}\\
			H_{A}:&\{\text{almost complex manifolds}\}\to\{\text{bigraded modules over }H_{BC}\}.
		\end{align*}
		together with natural transformations
		\[
		\begin{tikzcd}
			&H_{BC}\ar[ld]\ar[d]\ar[rd]&\\
			H_{Dol}\ar[rd]&(H_{dR},F,\bar F)\ar[d]&H_{Dolbar}\ar[ld]\\
			&H_A
		\end{tikzcd}
		\]
		and natural $\C$-antilinear isomorphisms $H_{BC}^{p,q}(X)\cong H_{BC}^{q,p}(X)$ and $H_{A}^{p,q}(X)\cong H_A^{q,p}(X)$.	
	\end{thm}
	\begin{rem}
		Fixing a hermitian metric, write $\cH_{BC}^{p,q}(X)=\cH_{BC}(X)\cap A_X^{p,q}$ and $\cH_A^{p,q}(X)=\cH_A(X)\cap A_X^{p,q}$ denote the pure-type harmonic forms in the sense of \cite{tardinitomassini20}. Nicoletta Tardini informed us that there are injections $\cH_{BC}^{p,q}(X)\subseteq H_{BC}^{p,q}(X)$ and $\cH^{p,q}_A(X)\subseteq H_A^{p,q}(X)$. 
	\end{rem}
	
	\begin{rem}[Infinite dimensionality]\label{rem: HBC infinite dimensional}
		If $X$ is a maximally non-integrable $2n$-manifold for $n=2,3$, then $H_{BC}^{n,n-1}(X)\cong H_{Dol}^{n,n-1}(X)$ and $H_{Dol}^{0,1}(X)\cong H_A^{0,1}(X)$. In particular, they may also be infinite dimensional on compact manifolds.
	\end{rem}
	
	\begin{rem}[a strong $\del\delbar$-property implies integrability]
		A complex manifold $X$ satisfies the so-called $\del\delbar$-property 
		if and only if all maps in the diagram in Proposition~\ref{prop: cohomology diagram} are isomorphisms. 
		One may use this as a definition of the $\del\delbar$-property for general almost complex structures. However, from Theorem~\ref{thm: mani => HDol infty} or Remark~\ref{rem: HBC infinite dimensional}, it follows that an almost complex structure on a $4$-manifold which satisfies the $\del\delbar$-property in this sense is already integrable.
	\end{rem}

\end{document}